\documentclass{amsart}

\usepackage[mathscr]{euscript}
\usepackage{epsfig}
\usepackage{latexsym}

\usepackage[all]{xy}

\usepackage{amssymb}
\usepackage{amsthm}
\usepackage{amsmath}
\usepackage{amsxtra}

.tfm
.tfm

\theoremstyle{plain}
\newtheorem{prop}{Proposition}[section]
\newtheorem{thm}[prop]{Theorem}
\newtheorem{coro}[prop]{Corollary}

\newtheorem{lemma}[prop]{Lemma}

\theoremstyle{definition}

\theoremstyle{remark}

\newtheorem*{remark}{Remark}

\numberwithin{table}{section}

\DeclareMathOperator{\Gal}{Gal}

\DeclareMathOperator{\lcm}{lcm}

\DeclareMathOperator{\coho}{H}

\newcommand{\F}{\mathbb F}

\newcommand{\Disc}{{\mathscr{D}}}

\newcommand{\Id}{id}
\newcommand{\Inf}{Inf}
\newcommand{\Res}{Res}

\def\NN{\mathbb N}
\def\RR{\mathbb R}

\def\<#1>{{\left\langle{#1}\right\rangle}}

\def\Z{{\mathbb Z}}             
\def\Q{{\mathbb Q}}             

\def\set#1{{\left\{{\def\st{\;:\;}#1}\right\}}}

\def\dkro#1#2{\left(\frac{#1}{#2}\right)}             
\let\kro\dkro

\def\O{{R}}           
\def\Ol(#1){{\mathop{\O_l}(\id{#1})}}                 
\def\Or(#1){{\mathop{\O_r}(\id{#1})}}                 
\def\Orp(#1){{\mathop{\O_r}(\idp{#1})}}               
\def\Otern(#1){{\mathop{\O^0_r}(\id{a})}}    

\def\id#1{{\mathfrak{#1}}}      
\def\idp#1{{\id{#1}_p}}         

\DeclareMathOperator{\trace}{{\mathrm{Tr}}}

\DeclareMathOperator{\M}{{\mathscr M}}

\let\MR\undefined
\DeclareMathOperator{\MR}{{\M_{\RR}}}

\def\T_#1(#2){{\mathop{\mathscr T}\nolimits_{#1}(\id{#2})}}
\def\TO(#1)_#2(#3){{\mathop{\mathscr T}\nolimits^{#1}_{#2}(\id{#3})}}

\def\A_#1(#2){{\mathop{\mathscr A}\nolimits_{#1}(\id{#2})}}
\def\Ax_#1(#2){{\mathop{\widetilde{\mathscr A}}\nolimits_{#1}(\id{#2})}}

\begin{document}
\title{On the embedding problem for $2^+S_4$ representations}
\author{Ariel Pacetti}
\address{Departamento de Matem\'atica, Universidad de Buenos Aires,
         Pabell\'on I, Ciudad Universitaria. C.P:1428, Buenos Aires, Argentina}
\email{apacetti@dm.uba.ar}
\thanks{The author was supported by a CONICET grant}
\thanks{The author would like to thank the "Universitat de Barcelona" where this work
was done}
\begin{abstract} Let $2^+S_4$ denote the double cover of $S_4$
  corresponding to the element in $\coho^2(S_4,\Z/2\Z)$ where
  transpositions lift to elements of order $2$ and the product of
  two disjoint transpositions to elements of order $4$ (denoted $\tilde
  S_4$ in \cite{Serre}). Given an elliptic curve $E$, let $E[2]$
  denote its $2$-torsion points. Under some conditions on $E$ (as in
  \cite{Bayer}) elements in $\coho^1(\Gal_\Q,E[2])\backslash \{ 0 \}$ correspond to
  Galois extensions $N$ of $\Q$ with Galois group (isomorphic to)
  $S_4$. On this work we give an interpretation of the addition law on
  such fields, and prove that the obstruction for $N$ having
  a Galois extension $\tilde N$ with $\Gal(\tilde N/ \Q) \simeq
  2^+S_4$ gives an homomorphism $s_4^+:\coho^1(\Gal_\Q,E[2])
  \rightarrow \coho^2(\Gal_\Q,\Z/2\Z)$. As a Corollary we can prove
  (if $E$ has conductor divisible by few primes and high rank) the
  existence of $1$-dimensional representations attached to $E$ and use
  them in some examples to construct $3/2$ modular forms mapping via
  the Shimura map to (the modular form attached to) $E$.
\end{abstract}

\maketitle

\section*{Introduction}
The study of modular forms of weight $1$ is equivalent to that of two
dimensional continuous faithful irreducible complex representations of
$\Gal_\Q := \Gal(\bar\Q/\Q)$ (see \cite{Deligne}). Looking at their
projectivization we have five different kinds: cyclic, dihedral,
$A_4$, $S_4$, $A_5$ (see \cite{Dornhoff}). If the image is cyclic then
the original representation is abelian, hence reducible. The dihedral
case corresponds to weight $1$ modular forms which are linear
combination of theta series attached to binary quadratic forms. The
``special ones'' are the last three cases. They are constructed using
different approaches (see \cite{JehanneA} for algorithms to construct
the $A_4$ and $A_5$ cases, and \cite{Jehanne} for the $S_4$ case). To
study the $S_4$ case, in \cite{Bayer} the next method is proposed: let
$E$ be an elliptic curve over $\Q$ with negative discriminant (if the
discriminant is positive the same method gives Maas forms), no
$2$-torsion points over $\Q$ and non-trivial Selmer $2$-group. The set
$\coho^1(\Gal_\Q,E[2]) \backslash \{0\}$ is in one to one
correspondence with fields $N$ with Galois group $S_4$ over $\Q$
containing $\Q(E[2])$. The obstruction for $N$ having a field
extension $\tilde N$ with Galois group over $\Q$ isomorphic to
$2^+S_4$ is an element in $\coho^2(\Gal_\Q,\Z/2\Z)$ (see \cite{Serre}
for a formula of the obstruction and \cite{Crespo},\cite{Crespo2} for
a method to compute a solution to the embedding problem when the
obstruction is trivial). This induces a map $s_4^+ :
\coho^1(\Gal_\Q,E[2]) \backslash \{0\} \rightarrow
\coho^2(\Gal_\Q,\Z/2\Z)$. The main result of this work is that if we
define $s_4^+(0) = 0$, then $s_4^+$ is a group homomorphism. As a
corollary all elliptic curves (in the above conditions) of conductor
$2^rp^s$ with $r,s \in \NN_0$ and $p$ a prime number such that the
$2$-Selmer group has rank at least two have a $1$-dimensional
representations with Galois group $2^+S_4$ attached to them. We end
this work with some examples of how using these weight $1$ modular
forms one can construct weight $3/2$ modular forms mapping via Shimura
(see \cite{Shimura}, Main Theorem) to the modular form (attached to)
$E$.

I would like to thank the University of Barcelona for its
hospitality and specially to Professor Pilar Bayer for useful
conversations and contributions to improve this work. I would also
like to thank Professor Jordi Quer for explaining me how he
computed the explicit solution to the embedding problems of the
examples in \cite{Bayer} and some other useful conversations as
well.

\section{Correspondence between $\coho^1(\Gal_\Q,E[2])$ and fields}

Let $E$ be an elliptic curve over $\Q$ with negative discriminant and
no $2$-rational points. The field $L = \Q(E[2])$ is a Galois extension
of $\Q$ with Galois group $S_3$.
Let $S$ be a finite set of primes containing $2$ and the primes
dividing the conductor of $E$ and denote $\coho^1(\Gal_\Q,E[2],S)$
the cocycles unramified outside $S$. Abusing notation we will
denote $\coho^1(\Gal_\Q,E[2],S)^\times$ the set
$\coho^1(\Gal_\Q,E[2],S) \backslash \{0\}$. By $\Disc(K)$ we
denote the discriminant of the field $K$.

\begin{prop} The elements in $\coho^1(\Gal_\Q,E[2],S)^\times$ are in
  one to one correspondence with fields $N$ such that $L \subset N$,
  $\Gal(N/\Q) = S_4$ and $N/\Q$ is unramified outside
  $S$. Furthermore, if $K \subset N$ is a degree $4$ extension of $\Q$
  whose normal closure is $N$ then $\Disc(K) = \Disc(L)$ in $\Q^\times
  / (\Q^\times )^2$.
\end{prop}

\begin{proof} The correspondence is Proposition 1.1 of \cite{Bayer}
(although they do not state the ramification condition). The main
idea (that we will need latter) is that if $\phi$ is a non-trivial
cocycle, then $\phi|_{\Gal_L}$ is a group homomorphism. If we
denote $N_\phi$ the fixed field of $\ker(\phi)$, then $L \subset
N_\phi$ and $\Gal(N_\phi/L) \simeq \Z/2\Z \oplus \Z/2\Z$ (the
Klein group). The ramification condition follows from this
isomorphism. For the last statement, note that since $S_4$ has a
unique normal subgroup of index $2$ (namely $A_4$), there is a
unique quadratic Galois subextension, namely
$\Q(\sqrt{\Disc(K)})$. Since $\Q(\sqrt{\Disc(L)})$ is another
quadratic Galois subextension they must be equal.
\end{proof}

\subsection{Field Addition Interpretation}

The group structure of $\coho^1(\Gal_\Q,E[2])^\times$ induces a
group structure on fields $N$ satisfying the above condition.
Using elementary field theory we will show a natural construction
of such addition. It is an easy group theory exercise to check
that $S_4 \simeq S_3 \ltimes (\Z/2\Z \oplus \Z/2 \Z)$ where, if we
denote $\{P_1,P_2,P_3\}$ the three nonzero elements of $\Z/2\Z
\oplus \Z/2 \Z$, the action is given by $\sigma . P_i =
P_{\sigma(i)}$.

\begin{lemma} Let $N_1,N_2$ be fields
corresponding to elements in $\coho^1(\Gal_\Q,E[2],S)^\times$,
then $N_1 \cap N_2 = L$.
\end{lemma}

\begin{proof} Clearly $L \subset N_1 \cap N_2$. Since $N_1 \cap N_2$ is
Galois over $\Q$, $\Gal(N_1/N_1\cap N_2)$ corresponds to a
non-zero normal subgroup of $S_4$ contained in $\Gal(N_1/L)$
(isomorphic to the Klein group), hence $\Gal(N_1/N_1\cap
N_2)\simeq \Gal(N_1/L)$.
\end{proof}

In particular $N_1 N_2$ is a Galois extension of $\Q$ with Galois
group of order $96$ and $\Gal(N_1N_2/L) \simeq \Gal(N_1/L)\oplus
\Gal(N_2/L) \simeq \Z/2\Z \oplus \Z/2\Z \oplus \Z/2\Z \oplus
\Z/2\Z$.

\begin{prop} $\Gal(N_1N_2/\Q) \simeq S_3 \ltimes
  \oplus_{i=1}^2 (\Z/2\Z \oplus \Z/2\Z)$, where if $u,v \in \Z/2\Z \oplus \Z/2\Z$
and $\sigma \in S_3$, $\sigma.(u,v) = (\sigma.u,\sigma.v)$.
\end{prop}

\begin{proof} Let $K_i \subset N_i$ be the fixed field of the subgroup $S_3$ looked as a
  subgroup of $S_4 = \Gal(N_i/\Q)$ (fixing the fourth element). We have the
diagram:
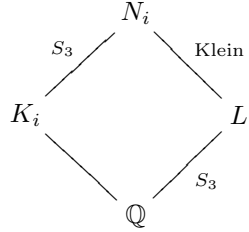
\[
\xymatrix{
& N_i \ar@{-}_{S_3}[ld] \ar@{-}[dr]^{\text{Klein}} &\\
K_i \ar@{-}[dr] & & L \ar@{-}^{S_3}[dl]\\
& \Q & }
\]

\noindent {\bf Claim:} $\Gal(N_1N_2/K_1K_2) \simeq \Gal(L/\Q)$. To
see this consider the sequence:
$$ \Gal(N_i/K_i) \hookrightarrow \Gal(N_i/\Q) \stackrel{\pi_i}{\rightarrow} \Gal(L/\Q)$$
where $\pi_i : \Gal(N_i/\Q) \rightarrow \Gal(L/\Q)$ is the
restriction map (which is the same as the quotient by
$\Gal(N_i/L)$). Since $\Gal(N_i/K_i)\cap \Gal(N_i/L) = \{\Id\}$,
the composition is an isomorphism.

Since $N_1K_2=N_1N_2$, $\Gal(N_1N_2/K_1K_2) \simeq \Gal(N_1/K_1)
\simeq \Gal(L/\Q)$ by the restriction map (respectively
$\Gal(N_1N_2/K_1K_2) \simeq \Gal(N_2/K_2)$). From the sequence
$$ \Gal(N_1N_2/K_1K_2) \hookrightarrow \Gal(N_1N_2/\Q) \stackrel{\Pi}{\rightarrow} \Gal(L/\Q)$$
given by restriction (where the composition is an isomorphism) we
conclude that $\Gal(N_1N_2/\Q) \simeq \Gal(L/\Q) \ltimes
\Gal(N_1N_2/L)$ and comparing the action with that of
$\Gal(N_i/\Q) \simeq \Gal(L/\Q) \ltimes \Gal(N_i/L)$, the result
follows.
\end{proof}

\begin{prop} The group $G:=S_3 \ltimes \sum_{i=1}^4 (\Z/2\Z)$ with the previous action has
three normal subgroups of order $4$. \label{subgroups}
\end{prop}

\begin{proof} Clearly the subgroups $\{0\} \ltimes (\Z/2\Z , \Z/2\Z
,0,0)$ , $\{0\} \ltimes (0,0,\Z/2\Z,\Z/2\Z)$ and $\{0\} \ltimes
\set{(a,b,a,b) \st a,b \in \Z/2\Z}$ are normal. Let $\coho
\triangleleft \,G$ be any normal subgroup of order $4$. Then $\Pi(H)
\triangleleft \, S_3$ with order $1$ or $2$. Since $S_3$ has no normal
subgroups of order $2$, $\Pi(H) = \{0\}$. The orbits of $S_3$ acting
on $\{0\} \ltimes \sum_{i=1}^4 (\Z/2\Z)$ are:
\begin{itemize}
\item $\set{ (0,0,0,0)}$
\item $\set{(1,0,0,0),(0,1,0,0),(1,1,0,0)}$
\item $\set{(0,0,1,0),(0,0,0,1),(0,0,1,1)}$
\item $\set{(1,0,1,0),(0,1,0,1),(1,1,1,1)}$
\item $\set{(1,0,0,1),(0,1,1,1),(1,1,1,0),(1,1,0,1),(1,0,1,1),(0,1,1,0)}$
\end{itemize}
\end{proof}

\begin{coro} Let $\psi_i \in \coho^1(\Gal_\Q,E[2],S)^\times$
and $N_i$ the corresponding field. The cocycle $\psi_1 + \psi_2$,
if non-trivial, corresponds to the field fixed by the third normal
subgroup of order $4$ in $\Gal(N_1N_2/\Q)$.
\end{coro}
\begin{proof}
The morphisms $\psi_i|_{\Gal_L}:\Gal_L \rightarrow \Z/2\Z \oplus
\Z/2\Z$, satisfy $N_i =\ker(\psi_i|_{\Gal_L})$. Clearly
$\psi_1+\psi_2$ is zero on $\Gal_{N_1N_2}$, hence its kernel is a
normal subgroup of order $4$ in $\Gal(N_1N_2/L)$ and the result
follows from Proposition \ref{subgroups}.
\end{proof}

\begin{remark} all normal subgroups of $G$ have pairwise trivial intersection
(corresponding to normal subfields $K_1$, $K_2$ and $K_3$). If we
define the subgroups:
\begin{itemize}
\item $\coho_4= \{0\} \ltimes \set{(0,0,0,0),(1,0,0,1),(0,1,1,1),(1,1,1,0)}$
\item $\coho_5= \{0\} \ltimes \set{(0,0,0,0),(1,1,0,1),(1,0,1,1),(0,1,1,0)}$
\end{itemize}
then all these five subgroups have trivial pairwise intersection.
They correspond to the other two (unique) subfields $N_4$ and
$N_5$ of index $4$ of $N_1N_2$ with the property that $N_i \cap
N_j = L$ for all $i \neq j$. Furthermore $N_4$ and $N_5$ are
Galois conjugates.

It is a nice exercise to prove that given any three order $4$
subgroups of $\oplus_{i=1}^4 \Z/2\Z$ having trivial pairwise
intersection, there exists another two order $4$ subgroups such
that all of them have the same property.
\end{remark}
\subsection{Two coverings of $S_4$}

We will consider cohomology groups with the trivial action. The
central $2$-extensions of $S_4$ correspond to elements in the
group $\coho^2(S_4,\Z/2\Z) \simeq \Z/2\Z \times \Z/2\Z$, where the
four groups are:
\begin{itemize}
\item $S_4 \oplus \Z / 2\Z$
\item $2^{\det}S_4$, corresponding to the cup product of the signature with itself.
\item $2^+S_4$ and $2^-S_4$.
\end{itemize}

The group $2^+S_4$ is isomorphic to $Gl_2(\F_3)$. A complete character
table of the group $2^+S_4$ can be found in \cite{Dornhoff}, Lemma
28.2.
\begin{lemma}
The group $2^+S_4$ has a subgroup isomorphic to $S_3$.
\label{keyingredient}
\end{lemma}
\begin{proof} The subgroup of $Gl_2(\F_3)$ spanned by $\< \left(\begin{array}{cc} 1&
 0\\ 0 & 1 \end{array} \right) , \left(\begin{array}{cc} 0& 1\\ 1 & 0
 \end{array} \right),  \left(\begin{array}{cc} 1& 2\\ 0 & 2
 \end{array} \right) >$ is isomorphic to $S_3$, an isomorphism given
  by sending the generators to the elements in $S_3$: $\Id$, $(1,2)$ and
$(1,3)$,
  respectively.
\end{proof}

Let $s_4^+$ denote the element in $\coho^2(S_4,\Z/2\Z)$
corresponding to the group $2^+S_4$.
Consider the projections from $G$ to $S_4$:
\begin{itemize}
\item $\Pi_1(\sigma,(x,y,z,w)) = (\sigma, (x,y))$.
\item $\Pi_2(\sigma,(x,y,z,w)) = (\sigma, (z,w))$.
\item $\Pi_3(\sigma,(x,y,z,w)) = (\sigma, (x+z,y+w))$.
\end{itemize}
where $\Pi_i$ maps $G$ to $\Gal(N_i/\Q)$ with kernel $H_i$, for
$i=1,2,3$ (the three normal subgroups of Proposition
\ref{subgroups}). The obstruction for the existence of $\tilde K_i$, a
field containing $K_i$ and Galois group $2^+S_4$ is the element in the
$2$-Brauer group $\Gamma^*(\Pi_i^*(s_4^+))$ where $\Gamma : \Gal_\Q
\rightarrow \Gal(K_1K_2/\Q)$ is the restriction map. Our main theorem can
be stated as

\begin{thm} $\Pi_1^*(s_4^+) + \Pi_2^*(s_4^+) + \Pi_3^*(s_4^+) =
  0$.
\label{main}
\end{thm}

From class field theory we know that the $2$-Brauer group injects
into the sum of its local components, i.e.
$\coho^2(\Gal_\Q,\Z/2\Z) \hookrightarrow \oplus_l
\coho^2(\Gal_{\Q_l},\Z/2\Z)$ and the local components are
 isomorphic to $\set{ \pm1}$. Let $\coho^2(\Gal_\Q,\Z/2\Z,S)$ be the subgroup of the
$2$-Brauer group of elements with trivial image at the primes
outside $S$. If we extend $s_4^+$ to $\coho^1(\Gal_\Q,E[2],S)$ by
setting $s_4^+(0)=0$, we get

\begin{coro} The map $s_4^+ : \coho^1(\Gal_\Q,E[2],S) \rightarrow
  \coho^2(\Gal_\Q,\Z/2\Z,S)$ is a group homomorphism.
\end{coro}
 From Serre's formula for the obstruction (\cite{Serre}, Theorem
1) it is clear that the image of $s_4^+$ is on this subgroup of
the $2$-Brauer group. The Corollary is an immediate consequence of
Theorem \ref{main} noting that the case when two fields are equal
is trivial from the fact that the cohomology groups are
$2$-groups.

\begin{proof}(Theorem)
Let $Z_2$ denote the group $\oplus_{i=1}^4 \Z/2\Z$ and consider the
  exact sequence
$$ 0 \rightarrow Z_2 \rightarrow G \rightarrow S_3 \rightarrow 0$$
using the inflation-restriction map we get an exact sequence
$$ 0 \rightarrow \coho^2(S_3,\Z/2\Z) \stackrel{\Inf}{\rightarrow}
\coho^2(G, \Z/2\Z) \stackrel{\Res}{\rightarrow} \coho^2(Z_2,\Z/2\Z) $$
Let $\psi := \Pi_1^*(s_4^+) + \Pi_2^*(s_4^+) +
\Pi_3^*(s_4^+)$. Restricted to $Z_2$, $\Pi_1 + \Pi_2 + \Pi_3 = 0$,
hence $\text{Res}(\psi) = 0$, i.e. $\psi$ is in the image of the
inflation map. This implies that $\psi$ does not depend on
representatives of the quotient map, in particular it is determined by
its values on $S_3 \ltimes \{(0,0,0,0)\} \times S_3 \ltimes
\{(0,0,0,0)\}$. By Lemma \ref{keyingredient}, $S_3$ is a subgroup of
$2^+S_4$, then $s_4^+((S_3,(0,0)),(S_3,(0,0))) = 0$ and the
result follows.
\end{proof}
\begin{remark} The fact that the group $2^+S_4$ has a subgroup
  isomorphic to $S_3$ is crucial for the map being a homomorphism. The
  same statement is false (in general) considering the maps between
  $\coho^1(\Gal_\Q,E[2],S) \rightarrow \coho^2(\Gal_\Q,\Z/2\Z,S)$ coming
  from the other groups $2^{\det}S_4$ and $2^-S_4$.
\end{remark}

Let $K_i$ be a degree $4$ extensions of $\Q$ with normal closure $N_i$
, $Q_{N_i}$ be the quadratic form $\trace_{K_i/\Q}(x^2)$ and
$W(Q_{N_i})$ its Witt invariant,

\begin{coro} $W(Q_{N_3}) = W(Q_{N_1}) + W(Q_{N_2}) +
  (2,\Disc(L))$ on $\text{Br}_2(\Q)$.
\end{coro}
\begin{proof}
This follows from Serre's Formula for the obstruction, see
\cite{Serre}, Theorem 1.
\end{proof}

\begin{coro} Let $E$ be an elliptic curve with conductor $2^rp^s$ with
  $r,s \in \NN_0$ and $2$-Selmer group of rank at least two, then
  there exists a $2^+S_4$ representation attached to $E$.
\label{existence}
\end{coro}
\begin{proof}
Since the Selmer group has rank at least $2$, let $N_1$ and $N_2$
be two different fields corresponding to elements in
 $\coho^1(\Gal_\Q,E[2])^\times$, and $N_3 := N_1+N_2$. Let $s_4^+(N_i)$
 denote their obstruction. From the injection of the $2$-Brauer group into its
local components it is clear that the obstruction is characterized
by the (finite set of) primes with $-1$ sign. Such set has an even
number of primes and is contained in the set $\set{2,p}$, hence if
two elements have non-zero obstruction in the Brauer group, the
third one does.
\end{proof}

\section{Applications and Examples}

We give a brief summary of how to construct the weight $1$ modular
forms attached to the Galois group $2^+S_4$. Let $K = \Q(x_1)$ be a
degree four extension of $\Q$ with normal closure $N$, an extension
with $S_4$ Galois group. By Theorem 1 of \cite{Serre} the obstruction
(for a lift with Galois group $2^+S_4$) is trivial if and only if the
quadratic form $\trace_{K/\Q}(x^2)$ is isomorphic (over $\Q$) to the
form $x_1^2+x_2^2+2x_3^2+2\Disc(K)x_4^2$. Furthermore if $P$ is a
transformation matrix sending one form to the other one, let
\[
\gamma := \det \left [
\left( \begin{array}{cccc}
1&x_1&x_1^2&x_1^3\\
1&x_2&x_2^2&x_2^3\\
1&x_3&x_3^2&x_3^3\\
1&x_4&x_4^2&x_4^3\\
\end{array} \right ) P \left( \begin{array}{cccc}
1 & 0 & 0 & 0\\
0 & 1 & 0 & 0\\
0 & 0 & \frac{1}{2} & \frac{1}{2}\\
0 & 0 & \frac{1}{2\sqrt{\Disc(K)}} & \frac{-1}{2\sqrt{\Disc(K)}}
\end{array} \right) + \Id \right ]
\]
If $\gamma \not = 0$, all the solutions to the embedding problem are
$\tilde N = N(\sqrt{c\gamma})$ with $c \in \Q^\times$ (see \cite{Crespo2},
Theorem 5). Furthermore, $c$ can be chosen such that $\tilde N$ is
unramified outside $S$ and has minimum ramification at the primes in
$S$. This choice of $c$ gives a weight $1$ modular, with character
$\kro{\Disc(K)}{}$ and minimum level. See \cite{Bayer}, Proposition
2.3, Proposition 2.4 and Proposition 2.5 to compute the Fourier
coefficients of the weight $1$ modular form once $\gamma$ is known.

We will give some examples of how this weight $1$ modular forms
can be used to construct some ``special'' $3/2$ modular forms.
Given a weight $2$ and level $p$ modular form $f$ (attached to an
elliptic curve $E$), in \cite{Gross} Gross gave a method to
construct a weight $3/2$ modular form (as linear combination of
theta series) in the Kohnen space mapping to $f$ via the Shimura
map. If $E$ has positive rank, the constructed weight $3/2$
modular form is the zero form. For this elliptic curves we will
show (in some examples) how using the weight $1$ modular form
coming from the solution of the obstruction problem (by Corollary
\ref{existence} we know that such form exists if $E$ has rank
greater than $1$) one can construct a non-zero weight $3/2$
modular form mapping to $f$ via the Shimura map. This construction
has some limitations (from a computation point of view) as we will
see latter, but works on many cases (our approach is similar to
that on \cite{Bungert}).

For $n \in \NN$, let
\[
\Theta_{n}(z) := \sum_{j=-\infty}^{+\infty} q^{nj^2}
\]
then the theta function $\Theta_n(z)$ is a weight $1/2$ modular form
of level $4n$ and character $\kro{n}{}$ (see \cite{Shimura} for the
definition of modular forms of half integral weight).

\begin{lemma}
Let $g(z)$ be a weight $1$, level $n$ and character $\kro{-d}{}$
modular form (with $d > 0$ and $d \mid n$), then $g(z) \Theta_d(z)$ is a modular
form of weight $3/2$, level $\lcm(n,4d)$ and trivial character.
\label{product}
\end{lemma}
\begin{proof} Let $M \in \Gamma_0(\lcm(n,4d))$, say $M = \left(
  \begin{array}{cccc}
\alpha & \beta\\
\gamma & \delta \end{array} \right)$, then:
\begin{itemize}
\item $g(Mz) = (\gamma z + \delta) \kro{-d}{\delta} g(z)$.
\item $\Theta_d(Mz)\Theta(z) = \kro{d}{\delta} \Theta_d(z) \Theta(Mz)$.
\item $\Theta(Mz)^2 = (\gamma z + \delta) \kro{-1}{\delta} \Theta(z)^2$.
\end{itemize}
Then $g(Mz) \Theta_d(Mz) \Theta(z)^3 = g(z) \Theta_d(z) \Theta(Mz)^3$
which is the definition of a weight $3/2$ modular form of trivial character.
\end{proof}

What we do is to compute the space of ternary quadratic forms
whose theta series are modular forms with trivial character (this
is equivalent to the forms have square discriminant), and add the
product of the weight $1$ modular form with the corresponding
theta series. Applying the Hecke operators on this set we look for
the eigenform corresponding to $E$. In practice the difference
between the dimension of the space of weight $3/2$ modular forms
and that of the subspace spanned by theta series increases with
the level (see \cite{Lehman} for some tables). Use this method in
general implies knowing too many Fourier coefficients of the
constructed modular forms, which is impracticable. For this reason
we use this method for the case $n=43$, where the dimension of the
whole space is almost the same as the subspace of theta series. If
the weight $1$ modular form $g$ has level $p$ and character
$\kro{-p}{}$, a prime congruent to $3 \, (\bmod 4$), $g(4z)$ has
level $4p$ and the same character, hence the product $g(4z)
\Theta_p(z)$ is in the Kohnen space of level $4p$ and trivial
character, which (by \cite{Kohnen} Theorem 1) is isomorphic to
$S_2(\Gamma_0(p))$. In these cases we can construct the weight
$3/2$ modular forms for big values of $p$.

\subsection{Examples} We use the method described above in some
particular examples. The computations were done with the PARI/GP
system~\cite{PARI}.

{\noindent \bf Notation:} the ternary quadratic forms will be denoted
$a_1$, $a_2$, $a_3$, $a_{23}$, $a_{13}$, $a_{12}$, to express the
form:
\begin{equation}\label{eq:notation:qf3}
  Q(X_1,X_2,X_3) = a_1 X_1^2 + a_2 X_2^2 + a_3 X_3^2 +
                   a_{23} X_2 X_3 + a_{13} X_1 X_3 + a_{12} X_1 X_2.
\end{equation}

\subsubsection{Case $43$} The elliptic curve $43A$ in Cremona's table
with equation $y^2+y=x^3+x^2$ has rank $1$. A generator is $P=(0,0)$
and it corresponds to the field $K$ with equation $P=x^4-2x-1$. By
\cite{Bayer} we know that the embedding problem for this case is
solvable, and a solution with minimal level is given by the element
\[
\gamma = 3(x_1^3x_2^2-x_2^2-x_1^2x_2+x_1x_2+x_2)+x_1^3-2x_1^2+4x_1
,
\]
where $x_1,x_2$ are roots of $P$. The corresponding modular form has
level $2^343$ and character $\kro{-43}{}$.
All ternary quadratic forms of level $2^243$ and $2^343$ with
trivial character are given in tables \ref{table172} and
\ref{table344}, where we denote $Q_{43\_i}$ (respectively
$Q_{86\_i}$) the $i$-th form on the table of ternary quadratic
forms of level $2^2 43$ (respectively of level $2^3 43$). See
\cite{Tornaria} for interactive tables of ternary quadratic forms
of a given level (in some specific genera) and \cite{Lehman},
Theorem $4$ and Theorem $5$, for the bijection between different
genera.

\begin{table}
\begin{tabular}{|l|rrrrrr||l|rrrrrr|}
\hline
 & $a_1$ & $a_2$ & $a_3$ & $a_{23}$ & $a_{13}$ & $a_{12}$ & &$a_1$ & $a_2$ & $a_3$ & $a_{23}$ & $a_{13}$ & $a_{12}$ \\
\hline\hline
$Q_1$ & $1$ & $11$ & $43$ & $0$ & $0$ & $1$ & $Q_{8}$ & $5$ & $18$ & $26$ & $18$ & $2$ & $4$\\
$Q_2$ & $4$ & $11$ & $14$ & $-10$ & $3$ & $2$ & $Q_{9}$ & $6$ & $15$ & $23$ & $2$ & $6$ & $4$\\
$Q_3$ & $6$ & $9$ & $10$ & $4$ & $5$ & $1$ & $Q_{10}$ & $9$ & $10$ & $24$ & $10$ & $2$ & $4$\\
$Q_4$ & $1$ & $43$ & $43$ & $0$ & $0$ & $0$ & $Q_{11}$ & $11$ & $14$ & $16$ & $-6$ & $4$ & $10$\\
$Q_5$ & $2$ & $22$ & $43$ & $0$ & $0$ & $2$ & $Q_{12}$ & $4$ & $43$ & $44$ & $0$ & $4$ & $0$\\
$Q_6$ & $3$ & $29$ & $29$ & $-28$ & $2$ & $2$ & $Q_{13}$ & $11$ & $16$ & $47$ & $16$ & $2$ & $4$\\
$Q_7$ & $4$ & $11$ & $43$ & $0$ & $0$ & $2$ & $Q_{14}$ & $15$ & $23$ & $24$ & $12$ & $8$ & $2$\\
\hline
\end{tabular}
\caption{Coefficients of ternary quadratic forms, level $2^2 43$.}
\label{table172}
\end{table}

\begin{table}
\begin{tabular}{|l|rrrrrr||l|rrrrrr|}
\hline
 & $a_1$ & $a_2$ & $a_3$ & $a_{23}$ & $a_{13}$ & $a_{12}$&  & $a_1$ & $a_2$ & $a_3$ & $a_{23}$ & $a_{13}$ & $a_{12}$ \\
\hline\hline
$Q_{1}$ & $3$ & $115$ & $115$ & $-114$ & $-2$ & $-2$ &$Q_{10}$ & $3$ & $29$ & $86$ & $0$ & $0$ & $-2$\\
$Q_{2}$ & $8$ & $43$ & $88$ & $0$ & $-8$ & $0$ & $Q_{11}$ & $5$ & $18$ & $86$ & $0$ & $0$ & $-4$\\
$Q_{3}$ & $19$ & $20$ & $91$ & $20$ & $6$ & $12$ & $Q_{12}$ & $8$ & $22$ & $43$ & $0$ & $0$ & $-4$\\
$Q_{4}$ & $4$ & $87$ & $87$ & $2$ & $4$ & $4$ & $Q_{13}$ & $19$ & $20$ & $26$ & $-4$ & $-16$ & $-12$\\
$Q_{5}$ & $15$ & $24$ & $92$ & $24$ & $4$ & $8$ & $Q_{14}$ & $1$ & $86$ & $86$ & $0$ & $0$ & $0$\\
$Q_{6}$ & $15$ & $23$ & $95$ & $-22$ & $-14$ & $-2$ & $Q_{15}$ & $6$ & $15$ & $86$ & $0$ & $0$ & $-4$\\
$Q_{7}$ & $16$ & $44$ & $47$ & $4$ & $16$ & $8$ & $Q_{16}$ & $9$ & $10$ & $86$ & $0$ & $0$ & $-4$\\
$Q_{8}$ & $23$ & $31$ & $47$ & $18$ & $14$ & $10$ & $Q_{17}$ & $13$ & $16$ & $47$ & $16$ & $6$ & $12$\\
$Q_{9}$ & $2$ & $43$ & $86$ & $0$ & $0$ & $0$ & $Q_{18}$ & $14$ & $21$ & $31$ & $14$ & $4$ & $12$\\
\hline
\end{tabular}
\caption{Coefficients of ternary quadratic forms, level $2^3 43$.}
\label{table344}
\end{table}

The theta functions of these quadratic forms are not linearly
independent. A basis is given by the theta functions of the
ternary quadratic forms: $\{Q_{43\_1}, Q_{43\_2}, \break
Q_{43\_3}, Q_{43\_4}, Q_{43\_5}, Q_{43\_6}, Q_{43\_7}, Q_{43\_8},
Q_{43\_9},Q_{43\_10},Q_{43\_11}, Q_{86\_1}, Q_{86\_2}, Q_{86\_3},
\break Q_{86\_4}, Q_{86\_5}, Q_{86\_6}, Q_{86\_7}, Q_{86\_8},
Q_{86\_9}, Q_{86\_11}\}$ (just by looking at enough Fourier
coefficients). The space $S_{3/2}(2^3 43)$ has dimension $25$ (see
\cite{Cohen} Theorem 2), hence there are $4$ modular forms
missing.

Let $f(z)$ denote the weight $1$ modular form of weight $1$, level
$2^3 43$ and character $\dkro{-43}{}$ associated to $K$. By Lemma
\ref{product}, $f \, \Theta_{43}$ is a weight $3/2$, level $2^3
43$ and trivial character cusp form. Since its coefficients are
not rational, we consider the two modular forms $F_1(z)=
\frac{1}{2}(f(z)+\bar{f}(z))\Theta_{43}$ and
$F_2(z)=\frac{\sqrt{-2}}{2}(f(z)-\bar{f}(z))\Theta_{43}$. These
two forms do have rational coefficients and are linearly
independent from the ternary theta functions. Looking at the $23$
modular forms together and computing the Hecke operators on them
we get the two eigenforms (we denote by $\Theta_{Q}$ the Theta
function of the ternary quadratic form $Q$):
\begin{enumerate}
\item $G_{43A} = -\frac{3}{2}\Theta_{Q_{43\_6}}+\frac{4}{3}\Theta_{Q_{43\_7}}+\Theta_{Q_{43\_8}}-\frac{16}{3}
\Theta_{Q_{86\_1}}+\Theta_{Q_{86\_2}}+6\Theta_{Q_{8\_3}}+\frac{5}{2}\Theta_{Q_{86\_9}}-5\Theta_{Q_{86\_11}}
-5F_2$.
\item $G_{172A} =-\frac{1}{2}\Theta_{Q_{43\_1}} -\frac{1}{2}\Theta_{Q_{43\_2}}+ \Theta_{Q_{43\_3}}+
\Theta_{Q_{43\_4}}-\Theta_{Q_{43\_5}}
+3\Theta_{Q_{43\_9}}+3\Theta_{Q_{43\_10}}
-6\Theta_{Q_{43\_11}}-2\Theta_{Q_{86\_4}}+ 4\Theta_{Q_{86\_5}}
+4\Theta_{Q_{86\_6}} -2\Theta_{Q_{86\_7}} -4\Theta_{Q_{86\_8}}
+\frac{3}{2}F_1$.
\end{enumerate}

They are Hecke eigenforms, and they map by the Shimura map to the
weight two modular forms associated to the elliptic curves $43A$ and
$172A$ on Cremona's table respectively. The curve $172A$ has rank $1$
and is given by the equation $y^2=x^3+x^2-13x+15$. The first
$50$ coefficients of their Fourier expansion are:\\

\noindent $\bullet$ $G_{43A}=q^2 + q^3 + q^5 - 5q^7 + 2q^8 - 4q^{12} - 3q^{18} + 2q^{19} + 2q^{20} -
3q^{22} - q^{26} - 2q^{27} + q^{29} + 2q^{32} - 3q^{33} + q^{34} + 4q^{37} +
5q^{39} + 2q^{42} + 2q^{43} - 3q^{45} + 3q^{46} + 6q^{48} -
3q^{50} + O(q^{51})$. \\

\noindent $\bullet$ $G_{172A}=q + q^6 - q^9 + q^{10} + q^{13} -
q^{14} - q^{17} + 2q^{21} - q^{25} + q^{41} - 3q^{49} +
O(q^{51})$.

\subsubsection{Case $563$:} the elliptic curve $563A$ with equation
$y^2+xy+y=x^3+x^2-15 x+16$ has rank $2$. The points $[2,-1]$ and
$[4,4]$ are generators for the rational points. The field
corresponding to the point $[2,-1]$ is given by the polynomial
$P=x^4 - 8 x^3 + 19 x^2 - 14 x - 1$. Its discriminant is $-563$
and the obstruction is trivial for this field. A solution to the
embedding problem is given by\\

$1126 \gamma=(57426x_2^2 - 408738x_2 - 155984)x_1^3 + (-434073 x_2^2 +
2329098 x_2 + 1542884) x_1^2 \break + (342834 x_2^2 - 1089141 x_2 - 4555297) x_1
-339522 x_2^3 + 2651994 x_2^2 - 6101295 x_2 + 4078271$,\\

\noindent where $x_1$ and $x_2$ are roots of $P$. Since $E$ has
discriminant $-563$, by Corollary 2.7 of \cite{Bayer} we know the
attached weight $1$ modular form has level $2^r563$ and character
$\kro{-563}{}$. By Theorem 2 of \cite{Rio} we know that $\gamma$
can be chosen such that the field $\tilde N$ above $N$ is
unramified at $2$ over $\Q$, hence the weight $1$ modular form has
level exactly $563$.

The field $\Q(\sqrt{\gamma})$ is given by $\Q(x_0)$ where $x_0$ is a
root of the polynomial:\\
\noindent $x^{24} - 3x^{23} - 9x^{22} + 22x^{21} + 55x^{20} -
68x^{19} - 212x^{18} + 85x^{17} + 467x^{16} - 34x^{15} - 698x^{14}
- 31x^{13} + 797x^{12} + 83x^{11} - 660x^{10} - 56x^9 + 420x^8 -
199x^6 + 32x^5 + 55x^4 - 20x^3 - 4x^2 + 3x + 1$.

Its discriminant is $-563^{11}$ (confirming that our choice of
$\gamma$ gives an extension unramified at $2$). All the Fourier
coefficients of this weight $1$ modular form can be computed as
stated before except the one corresponding to the ramified prime.
To compute $a_{563}$ we look at the inertia degree of $563$ in
$\Q(x_0)$ and since it is $2$ it follows that $a_{563}=-1$.

Let $F_{563}(z)$ denote the weight $1$ and level $563$ modular form
attached to this representation. The form $F_{563}(4z)
\Theta_{563}(z)$ is in the Kohnen space with trivial character, whose
space has dimension $48$ (while the whole space has dimension
$143$). The form $f_{563}(z) := \frac{1}{2}(F_{563}(4z)+\overline
F_{563}(4z)) \Theta_{563}(z)$, has rational coefficients. The
space spanned by the $33$ modular forms obtained from ternary
quadratic forms with trivial character in the Kohnen space (see
table \ref{table563}) and $f_{563}(z)$ is closed for the Hecke
operators, and the form\\

$F_{563A}=-11\Theta_{Q_1} -2\Theta_{Q_2}+ 44\Theta_{Q_3}
-8\Theta_{Q_4} -30\Theta_{Q_5}+ 16\Theta_{Q_6} -38\Theta_{Q_7}
+34\Theta_{Q_8} -2\Theta_{Q_9} +15\Theta_{Q_{10}}
+2\Theta_{Q_{11}} +22\Theta_{Q_{12}}-22\Theta_{Q_{14}}
+4\Theta_{Q_{15}} -7\Theta_{Q_{16}} -18\Theta_{Q_{17}}+
29\Theta_{Q_{18}} +26\Theta_{Q_{19}} +24\Theta_{Q_{20}}
+14\Theta_{Q_{21}} -6\Theta_{Q_{22}} +28\Theta_{Q_{23}}
-28\Theta_{Q_{24}} -34\Theta_{Q_{25}} -46\Theta_{Q_{26}}
+22\Theta_{Q_{27}} +8\Theta_{Q_{28}} -14\Theta_{Q_{29}}+
20\Theta_{Q_{30}}
-42\Theta_{Q_{31}}+13f_{563}$\\

\noindent is an eigenform for the Hecke operators mapping via Shimura
to the modular form (attached to the elliptic curve) $563A$. The first
$50$ coefficients of its Fourier expansion are:\\

$-2q^3 + 2q^4 - 2q^7 + 2q^{11} - 2q^{16} + 4q^{23} + 2q^{27} +
4q^{28} + 2q^{39} - 2q^{40} + 2q^{47} + 4q^{48} + O(q^{51})$.

\begin{table}
{\small
\begin{tabular}{|l|rrrrrr||l|rrrrrr|}
\hline
 & $a_1$ & $a_2$ & $a_3$ & $a_{23}$ & $a_{13}$ & $a_{12}$&  & $a_1$ & $a_2$ & $a_3$ & $a_{23}$ & $a_{13}$ & $a_{12}$ \\
\hline\hline

$Q_1$ & $4$ & $563$ & $564$ & $0$ & $-4$ & $0$ &$Q_{18}$ & $36$ & $68$ & $563$ & $0$ & $0$ & $-28$\\
$Q_2$ & $3$ & $751$ & $751$ & $-750$ & $-2$ & $-2$ & $Q_{19}$ & $23$ & $196$ & $299$ & $-96$ & $-22$ & $-4$\\
$Q_3$ & $39$ & $59$ & $584$ & $-52$ & $-32$ & $-14$& $Q_{20}$ & $40$ & $119$ & $284$ & $8$ & $20$ & $32$\\
$Q_4$ & $39$ & $67$ & $580$ & $-48$ & $-20$ & $-38$& $Q_{21}$ & $47$ & $100$ & $296$ & $-36$ & $-40$ & $-28$\\
$Q_5$ & $44$ & $52$ & $563$ & $0$ & $0$ & $-12$ & $Q_{22}$ & $68$ & $71$ & $299$ & $-62$ & $-16$ & $-36$\\
$Q_6$ & $47$ & $48$ & $575$ & $48$ & $2$ & $4$ & $Q_{23}$ & $44$ & $155$ & $207$ & $106$ & $20$ & $16$\\
$Q_7$ & $48$ & $51$ & $575$ & $14$ & $48$ & $28$ & $Q_{24}$ & $47$ & $155$ & $192$ & $-44$ & $-8$ & $-46$\\
$Q_8$ & $51$ & $52$ & $576$ & $52$ & $20$ & $40$ & $Q_{25}$ & $39$ & $179$ & $231$ & $174$ & $2$ & $30$\\
$Q_9$ & $7$ & $323$ & $644$ & $-320$ & $-4$ & $-6$ & $Q_{26}$ & $59$ & $120$ & $191$ & $40$ & $6$ & $36$\\
$Q_{10}$ & $12$ & $188$ & $563$ & $0$ & $0$ & $-4$ & $Q_{27}$ & $71$ & $107$ & $191$ & $-26$ & $-14$ & $-58$\\
$Q_{11}$ & $11$ & $207$ & $615$ & $-202$ & $-6$ & $-10$ & $Q_{28}$ & $75$ & $92$ & $215$ & $-84$ & $-38$ & $-24$\\
$Q_{12}$ & $16$ & $143$ & $567$ & $6$ & $16$ & $12$ & $Q_{29}$ & $63$ & $143$ & $144$ & $36$ & $16$ & $2$\\
$Q_{13}$ & $19$ & $119$ & $596$ & $-116$ & $-16$ & $-6$ & $Q_{30}$ & $71$ & $127$ & $160$ & $96$ & $20$ & $6$\\
$Q_{14}$ & $23$ & $99$ & $591$ & $-94$ & $-18$ & $-10$ & $Q_{31}$ & $76$ & $119$ & $160$ & $64$ & $60$ & $12$\\
$Q_{15}$ & $27$ & $84$ & $584$ & $84$ & $4$ & $8$ & $Q_{32}$ & $64$ & $143$ & $176$ & $-140$ & $-4$ & $-24$\\
$Q_{16}$ & $28$ & $84$ & $563$ & $0$ & $0$ & $-20$ & $Q_{33}$ & $103$ & $108$ & $171$ & $104$ & $86$ & $92$\\
$Q_{17}$ & $36$ & $63$ & $572$ & $4$ & $36$ & $8$ & & & & & & &\\
\hline
\end{tabular}
\caption{Coefficients of ternary quadratic forms, level $2^2 563$.}
\label{table563}
}
\end{table}

\subsubsection{Case $643$} This case is similar to the
previous one. The elliptic curve $643A$ given by the equation
$y^2+xy=x^3-4x+3$ has rank $2$ and the points $[1,0]$ and $[2,1]$
generate the rational points. Their obstruction is non-trivial.
They sum is the point $[-1,3]$ which do have trivial obstruction.
It corresponds to the polynomial $P = x^4 - x^3 - 2x + 1$, whose
field has discriminant $-643$. A solution to the embedding problem
is given by the element:\\

$643 \gamma = (-123456 x_2^2 + 36008 x_2 - 1376) x_1^3 + (79732 x_2^2
- 70820 x_2 + 21952) x_1^2 + (-9092 x_2^2 - 101504 x_2 + 51440) x_1
-75964 x_2^3 + 176272 x_2^2 + 43724 x_2 - 25540$,\\

\noindent where $x_1$ and $x_2$ are roots of $P$. The field
$\Q(\sqrt{\gamma})$ is given by $\Q(x_0)$ where $x_0$ is a root of the
polynomial:\\
\noindent $x^{24} - 5 x^{23} + 11 x^{22} - 8 x^{21} - 10 x^{20} +
23 x^{19} + 9 x^{18} - 86 x^{17} + 171 x^{16} - 121 x^{15} - 212
x^{14} + 636x^{13} - 504 x^{12} - 156 x^{11} + 766 x^{10} - 1116
x^9 + 1364 x^8 - 1100 x^7 + 697 x^6 - 426 x^5 + 227 x^4 - 37 x^3 +
25 x^2 - 29 x + 5$.

 This field has discriminant $-643^{11}$ hence the weight $1$ modular
form has level $643$ and character $\kro{-643}{}$. Looking at the
inertia, its Fourier coefficient $a_{643} = -1$.

Let $F_{643}(z)$ denote the weight $1$ and level $643$ modular form
attached to this representation. The form $F_{643}(4z)
\Theta_{643}(z)$ is in the Kohnen space with trivial character, whose
space has dimension $54$ (while the whole space has dimension
$163$). The form $f_{643}(z) := \frac{1}{2}(F_{643}(4z)+\overline
F_{643}(4z)) \Theta_{643}(z)$, has rational coefficients. The
space spanned by the $30$ modular forms obtained from ternary
quadratic forms with trivial character in the Kohnen space (see
table \ref{table643}) and $f_{643}(z)$ is closed for the Hecke
operators, and the form\\

$F_{643A} = -3\Theta_{Q_1} -\Theta_{Q_2}+ 6\Theta_{Q_3}
-2\Theta_{Q_4} -\Theta_{Q_5} +8\Theta_{Q_6} +\Theta_{Q_7}+
7\Theta_{Q_8}+ 3\Theta_{Q_9} -7\Theta_{Q_{10}} +\Theta_{Q_{11}}
-4\Theta_{Q_{12}} +4\Theta_{Q_{13}} +6\Theta_{Q_{14}}
+12\Theta_{Q_{15}} -11\Theta_{Q_{16}}+ 6\Theta_{Q_{17}}
-4\Theta_{Q_{18}} -3\Theta_{Q_{19}} -\Theta_{Q_{20}}
-9\Theta_{Q_{21}} +2\Theta_{Q_{22}} +2\Theta_{Q_{24}}
-8\Theta_{Q_{25}}
-12\Theta_{Q_{27}} +2\Theta_{Q_{28}} +6\Theta_{Q_{29}}+ 4f_{643}$\\

\noindent is an eigenform for the Hecke operators mapping via Shimura
to the modular form (attached to the elliptic curve) $643A$. The first
$50$ coefficients of its Fourier expansion are:\\

$q^4 - q^7 + q^{15} - q^{16} + q^{23} - q^{24} + 2q^{28} + q^{31}
- q^{36} + q^{40} + O(q^{51})$.

\begin{table}
{\small
\begin{tabular}{|l|rrrrrr||l|rrrrrr|}
\hline
 & $a_1$ & $a_2$ & $a_3$ & $a_{23}$ & $a_{13}$ & $a_{12}$&  & $a_1$ & $a_2$ & $a_3$ & $a_{23}$ & $a_{13}$ & $a_{12}$ \\
\hline\hline
$Q_1$ & $4$ & $643$ & $644$ & $0$ & $-4$ & $0$ & $Q_{16}$ & $55$ & $104$ & $328$ & $52$ & $12$ & $48$ \\
$Q_2$ & $7$ & $368$ & $735$ & $368$ & $2$ & $4$ & $Q_{17}$ & $63$ & $88$ & $327$ & $44$ & $10$ & $40$ \\
$Q_3$ & $16$ & $163$ & $647$ & $6$ & $16$ & $12$ & $Q_{18}$ & $63$ & $95$ & $335$ & $62$ & $46$ & $58$ \\
$Q_4$ & $23$ & $112$ & $671$ & $112$ & $2$ & $4$ & $Q_{19}$ & $60$ & $131$ & $216$ & $4$ & $20$ & $24$ \\
$Q_5$ & $28$ & $92$ & $643$ & $0$ & $0$ & $-4$ & $Q_{20}$ & $51$ & $152$ & $256$ & $-148$ & $-28$ & $-12$ \\
$Q_6$ & $31$ & $83$ & $671$ & $-82$ & $-30$ & $-2$ & $Q_{21}$ & $83$ & $95$ & $220$ & $36$ & $32$ & $26$ \\
$Q_7$ & $15$ & $343$ & $344$ & $172$ & $8$ & $2$ & $Q_{22}$ & $92$ & $95$ & $231$ & $74$ & $52$ & $64$ \\
$Q_8$ & $24$ & $215$ & $323$ & $2$ & $12$ & $8$ & $Q_{23}$ & $64$ & $167$ & $204$ & $-152$ & $-28$ & $-40$ \\
$Q_9$ & $23$ & $228$ & $339$ & $-104$ & $-18$ & $-20$ & $Q_{24}$ & $96$ & $136$ & $139$ & $60$ & $44$ & $28$ \\
$Q_{10}$ & $31$ & $168$ & $332$ & $84$ & $4$ & $16$ & $Q_{25}$ & $96$ & $135$ & $156$ & $-44$ & $-92$ & $-20$ \\
$Q_{11}$ & $40$ & $135$ & $324$ & $8$ & $20$ & $32$ & $Q_{26}$ & $111$ & $116$ & $143$ & $-68$ & $-42$ & $-8$ \\
$Q_{12}$ & $39$ & $132$ & $339$ & $-64$ & $-38$ & $-4$ & $Q_{27}$ & $116$ & $136$ & $143$ & $76$ & $68$ & $108$ \\
$Q_{13}$ & $39$ & $135$ & $331$ & $-62$ & $-14$ & $-22$ & $Q_{28}$ & $104$ & $124$ & $167$ & $104$ & $88$ & $12$ \\
$Q_{14}$ & $36$ & $143$ & $359$ & $-142$ & $-16$ & $-4$ & $Q_{29}$ & $111$ & $131$ & $143$ & $-54$ & $-42$ & $-82$ \\
$Q_{15}$ & $55$ & $95$ & $332$ & $52$ & $32$ & $18$ & $Q_{30}$ & $104$ & $144$ & $167$ & $-60$ & $-88$ & $-92$ \\
\hline
\end{tabular}
\caption{Coefficients of ternary quadratic forms, level $2^2 643$.}
\label{table643}
}
\end{table}

\def\MR#1{}
\bibliography{bibliography}
\bibliographystyle{amsplain}

\end{document}